\documentclass{amsart}
\usepackage{amssymb}
\DeclareMathSymbol{\twoheadrightarrow} 
{\mathrel}{AMSa}{"10}

\def\Q{{\mathbf Q}}
\def\Z{{\mathbf Z}}

\def\Gal{\mathrm{Gal}}
\def\End{\mathrm{End}}
\def\Aut{\mathrm{Aut}}

\def\Hom{\mathrm{Hom}}

\def\GL{\mathrm{GL}}

\def\Sp{\mathrm{Sp}}

\def\dim{\mathrm{dim}}

\def\O{{\mathcal O}}

\def\bmu{\boldsymbol \mu}
\newtheorem{thm}{Theorem}[section]
\newtheorem{lem}[thm]{Lemma}
\newtheorem{cor}[thm]{Corollary}
\newtheorem{prop}[thm]{Proposition}
\theoremstyle{definition}

\title[Symplectic representations of inertia groups]
{Symplectic representations of inertia groups}

\thanks{Silverberg would like to thank NSA and the
Alexander-von-Humboldt Stiftung for financial support,
and Herbert Lange and 
the Mathematics Institute of the University of Erlangen 
for their hospitality. 
Zarhin would like to thank NSF for financial support.}

\author[A.\ Silverberg]{A.\ Silverberg}
\address{Mathematics Department, 
Ohio State University,
231 W.\ 18 Avenue,
Columbus, Ohio, 43210,
USA}
\email{silver\char`\@math.ohio-state.edu}

\author[Yu. G. Zarhin]{Yu. G. Zarhin}
\address{Mathematics Department, 
Pennsylvania State University,
University Park, PA, 16802,
USA\newline\indent
Institute for Mathematical Problems in Biology,
Russian
Academy of Sciences, Push\-chi\-no, Moscow Region,
142292,
Russia}
\email{zarhin\char`\@math.psu.edu}

\begin{document}

\maketitle

\section{Introduction and notation}

The aim of this paper is to study finite 
inertia subgroups of symplectic groups over the field 
 $\Q_\ell$ of $\ell$-adic numbers. A finite group 
is called an inertia group with respect to 
a given prime $p \ne \ell$ if it is a semi-direct product of 
a finite normal $p$-subgroup and a cyclic $p'$-group. 
These groups are exactly the inertia groups of finite 
Galois extensions of discrete valuation fields of residue 
characteristic $p$. The study of semistable reduction of 
abelian varieties over such fields leads naturally to certain 
finite inertia subgroups of the
symplectic group $\Sp_{2g}(\Q_\ell)$ \cite{SZAl}. 
In \cite{SZcomp} we constructed examples for every odd prime 
$\ell$ of inertia subgroups of $\Sp_{2g}(\Q_\ell)$ which are 
not conjugate, even in 
$\GL_{2g}(\Q_\ell)$, 
to a subgroup of $\Sp_{2g}(\Z_\ell)$.  However, 
it turns out (and this is the main result of this paper) 
that every finite inertia subgroup of $\Sp_{2g}(\Q_\ell)$ is 
{\em isomorphic} to a subgroup of $\Sp_{2g}(\Z_\ell)$, if $\ell>3$. 

See \cite{Fontaine} for a study of representations of 
inertia groups in characteristic $0$.

Throughout this paper $\ell$ is an odd prime,
$K$ is a field that is an unramified finite  
extension of $\Q_\ell$, and  
$G$ is a finite group that is a semi-direct product
	$G=H L$
of a normal $\ell'$-subgroup $H$ and a cyclic $\ell$-group $L$. 
Note that if $G$ is a finite inertia group for some prime $p$, 
then $G$ is of this form for every prime $\ell \ne p$ (see 
Lemmas 3.2 and 3.3 of \cite{SZcomp}).
We always assume that 
the group algebra $K[H]$ is decomposable, i.e.,
splits into a direct sum of matrix algebras over fields
 (this will automatically
be the case when $G$ is an inertia group; see \cite{SerreArtin}).
We write $\Sp_{2d}(R)$ for the group of $2d \times 2d$
symplectic matrices over a ring $R$.
If $E$ is a field that is a finite extension of $\Q_\ell$, let
$\O_E$ denote the ring of integers.

The following result is the main result of this paper.

\begin{thm}[Embedding Theorem]
Suppose $d$ is a positive integer, 
and there is an embedding
$G \hookrightarrow \Sp_{2d}(K)$.
If $\ell \ge 5$, then 
there is an embedding
	$G \hookrightarrow \Sp_{2d}(\O_K)$.
\end{thm}

We will make use of the following result, 
which we prove in \S\ref{Hproof}.
Our proof was inspired by \S 17.6 
in \cite{SerreRep} and the proof of Lemma 1.1 in
\S 1 of Chap.~X in \cite{Feitbook}.

\begin{thm}[Extension Theorem]
\label{Htheorem}
Suppose that $W$ is a finite dimensional $K$-vector space,
$f: W \times W \to K$
is a non-degenerate alternating (resp., symmetric) 
$K$-bilinear form, and
$\tau: H \to \Aut(W,f)$
is a group homomorphism that makes $W$ into a simple $K[H]$-module. 
Assume that for every $g \in G$, the representation of $H$
$$\tau_g: H \to \Aut(W), \quad h \mapsto \tau(g h g^{-1})$$
is isomorphic to $\tau$.
Suppose $T$ is an $H$-stable $\O_K$-lattice in $W$. 
Then  
$$\tau: H \to \Aut(T,f) \subset \Aut(W,f)$$
can be extended to a homomorphism
$$\tau_G: G \to \Aut(T,f) \subset \Aut(W,f).$$
\end{thm}

We will say that two bilinear forms have the same parity if
they are either both symmetric or both alternating.
Let $\pi$ denote a uniformizer of $\O_K$.

\section{Proof of the Extension Theorem}
\label{Hproof}

Replacing $f$ by $\pi^i f$ if necessary, 
we may assume that $f(T,T)=\O_K$ and
$$f: T \times T \to \O_K$$
is perfect.

Let $E$ denote the centralizer of $H$ in $\End_K(W)$.
Since $K[H]$ is decomposable and $W$ is a simple
$K[H]$-module, $E$ is a field that is a finite 
extension of $K$.
It follows from Theorem 74.5 
(especially a description of the center K in (ii)) 
of \cite{CR2} that
$E$ is the field of definition of a certain character of $H$; 
in particular,
	$E \subseteq K(\bmu_{\# H})$.
Thus $E/K$ is unramified, since $\# H$ is prime to $\ell$.
Therefore $E/\Q_\ell$ is also unramified, since $K/\Q_\ell$ is unramified.

Since the $H$-module W is simple and $H$ is an $\ell'$-group,
the $H$-module $T/\pi T$ is also simple. 
This implies that the $H$-stable
$\O_K$-sublattice $\O_E T$ of $W$ is of the form $aT$ with $a \in K^*$ 
(see Exercise 15.2 
of \cite{SerreRep}). 
Thus $T= a^{-1}(\O_E T)$ is $\O_E$-stable, 
i.e., $T$ is an $\O_E$-lattice in the $E$-vector space $V$.

It follows from the Jacobson density theorem that the image of
$K[H]$ in $\End_K(W)$ is $\End_E(W)$. In other words,
$\End_E(W)$ is the $K$-vector subspace of $\End_K(W)$
generated by the $\tau(h)$ for $h \in H$.

By the non-degeneracy of $f$, there is an involution
$u \mapsto u'$ of $\End_K(W)$ characterized by 
$$f(ux,y)=f(x,u'y) \quad \forall x,y \in W.$$
By the $H$-invariance of $f$,
$\tau(h)'=\tau(h^{-1}) \quad \forall h \in H$.
Thus the involution $u\mapsto u'$ sends $\End_E(W)$ into itself, 
and therefore sends $E$, the center of $\End_E(W)$, 
into itself.
Let 
$$E_0=\{u\in E \mid u'=u\}.$$
Then 
$K \subseteq E_0 \subseteq E$.
Either $E=E_0$ or $E/E_0$ is a quadratic extension.

Let $c$ be a generator of $L$. 
The homomorphisms
$$\tau: H \to \Aut_{\O_K}(T), \quad \tau_c: H \to \Aut_{\O_K}(T)$$
define $\O_K[H]$-module structures on $T$
in such a way that the corresponding
$K[H]$-modules are isomorphic. Since $H$ is an $\ell'$-group, 
the corresponding $\O_K[H]$-modules are isomorphic
(see \S 14.4 and \S 15.5 of \cite{SerreRep}),
i.e., there exists
$A \in \Aut_{\O_K}(T)$ such that
$$\tau(c h c^{-1})=A \tau(h) A^{-1} \quad \forall h \in H.$$
Then
$$\tau(c^i h c^{-i})=A^i \tau(h) A^{-i} \quad \forall h \in H, i \in \Z.$$
Since $c^{\# L}=1$, we have $A^{\# L} \in \O_E^*$. Further, 
$ A\End_E(W)A^{-1}= \End_E(W)$.
Since $E$ is the center of $\End_E(W)$,
$A E A^{-1} =E$.
In other words, defining 
$\iota(c)(u)=AuA^{-1}$ for $u \in E$ induces a homomorphism
	$$\iota: L \to \Gal(E/K).$$
Since $[E:E_0]$ divides $2$ and 
$\ell$ is odd, the composition
$L \to \Gal(E/K) \twoheadrightarrow \Gal(E_0/K)$
has the same kernel as $\iota$.
Thus 
\begin{equation}
\label{idiv}
\#\iota(L) \text{ divides } [E_0:K].
\end{equation}

Let 
$f_A(x,y)= f(Ax,Ay)$ for $x, y \in T$.
 For all $h \in H$,
$$f_A(\tau(h)x,\tau(h)y)=f(A\tau(h)x, A\tau(h)y)=$$
$$f(A \tau(h) A^{-1} Ax,A \tau(h)A^{-1} Ay)=
f(\tau(c h c^{-1}) Ax, \tau(c h c^{-1}) Ay)=f_A(x,y),$$
since $c h c ^{-1} \in H$ and $f$ is
$H$-invariant.  
Thus $f_A$ is $H$-invariant and of the same parity as $f$.  
Therefore there exists $a \in E_0^*$ such that
$$	f(Ax,Ay)= f(ax,y)=f(x,ay) \quad \forall x,y \in W.$$
Since $A\in\Aut(T)$ and $f:T\times T\to \O_K$ is perfect, we have
$a\in\Aut(T)$. This implies easily that 
$a \in \O_{E_0}^*$.

Let $\sigma := \iota(c) \in \Gal(E/K)$.
Then $aA=A\sigma^{-1}(a)$.
There exists $a_1 \in \O_{E_0}^* $ such that
	$$a \sigma^{-1}(a)=a_1^2.$$
(Indeed, let $\eta$ be a uniformizer for $\O_{E_0}$.
Then $\sigma(a)^{-1} \equiv a^{\ell^j} \pmod{ \eta}$ 
for some non-negative integer $j$. 
Since $\ell$ is odd, $\ell^j +1$ is even, so 
$a\sigma^{-1}(a)$ is a square modulo $ \eta$.
Thus $a\sigma^{-1}(a)$ is a square, since all elements of 
$\O_{E_0}$ congruent to $1$ modulo $ \eta$ are squares.)

For all $x,y \in T$,
$$f((A^2 a_1^{-1})x,(A^2 a_1^{-1})y)=
f(aAa_1^{-1}x,A a_1^{-1}y)=$$
$$f(A\sigma^{-1}(a)a_1^{-1} x, A a_1^{-1}y)=
f(a \sigma^{-1}(a)a_1^{-1}x, a_1^{-1}y)= f(x,y).$$
Let $A_1 = A^2a_1^{-1} \in \Aut_{\O_K}(T)$. 
Then $f$ is $A_1$-invariant,
$\det(A_1)=1$, 
and conjugation by $A_1$  coincides with conjugation by $c^2$. 
Thus ${A_1}^{\# L}= bI$ is a scalar operator
 of determinant 1 on the $E_0$-vector space $W$. Therefore
$b \in E_0^*$ is a root of unity, i.e., 
${b}^\mu=1$ where $\mu$ is the number of
roots of unity in $E_0$. 
Since $E_0 \subseteq E$, the extension $E_0/\Q_\ell$
is unramified, so $\ell$ does not divide $\mu$.
Letting
	$B:={A_1}^\mu$,
then conjugation by B coincides with conjugation
by $b_1:=c^{2\mu}$, and $B^{\# L}=I$. Since $b_1$ is a generator
of $L$, sending $b_1$ to $B$ defines the desired extension
$\tau_G$ of $\tau$. 

\section{Lemmas for the Embedding Theorem}

Assume from now on that we are in the setting of the Embedding Theorem.
Therefore there exist a $2d$-dimensional $K$-vector space $V$,
a non-degenerate alternating $K$-bilinear form 
$$e: V \times V \to K,$$
and a faithful symplectic representation
$$\rho: G \hookrightarrow \Aut(V,e).$$

\begin{prop}
\label{isotype}
Suppose that $V$ is simple as a G-module but not as an $H$-module.
Then either 
\begin{enumerate}
\item[(i)] the $H$-module $V$ is isomorphic
to $W^r$ for some simple $H$-module $W$ and some $r>1$, 
or 
\item[(ii)] there exist a normal subgroup $G_1$
of $G$, and a simple symplectic $G_1$-module $V_1$
which is a $K$-vector space of dimension $2d/[G:G_1]$, 
such that $H \subseteq G_1 \ne G$
and such that if $g_1$ is a non-identity element of $G_1$,
then there exists $g \in G$ such that $g g_1 g^{-1}$ is not
in the kernel of $G_1 \to \Aut(V_1)$.
\end{enumerate}
\end{prop}

\begin{proof}
We follow the proof of Prop.\ 24 
of \cite{SerreRep}. 
Let $V=\oplus_{i=1}^n V_i$ be the canonical decomposition of the 
restriction of $\rho$ to $H$  
into a direct sum of isotypic representations. 
Since the $G$-module $V$ is simple, $G$ permutes the $V_i$ transitively. 
If $V$ is some $V_i$, then the $H$-module $V$ 
is isotypic and (i) holds.
Assume from now on that (i) does not hold. Let
$$G_1=\{s\in G\mid s(V_1)=V_1\}\subset G.$$
Then $H \subseteq G_1\ne G$. 
Since $G_1$ contains $H$, it is normal in $G$. 
Thus for every $V_i$,
$$G_1=\{s\in G\mid s(V_i)=V_i\}\subset G.$$
Every $V_i$ is a simple $G_1$-module, 
because $V$ is a simple $G$-module. 

The kernels of the natural maps $G_1 \to \Aut(V_i)$
have trivial intersection and are conjugate in $G$.
It follows that if $g_1$ is a non-identity element of $G_1$,
then it has a conjugate which does not lie 
in the kernel of $G_1 \to \Aut(V_1)$.
  
Since $[G:G_1]$ divides $[G:H]$, it 
is an $\ell$-power, and therefore odd. Thus
$n$ is odd, so at least one of the 
simple $G_1$-modules $V_i$ is self-dual.
This implies easily that all the $V_i$ are self-dual.
Suppose $V_1$ is not symplectic. Then none of the $V_i$
are symplectic. Since $V_i$ is a simple 
$G_1$-module, every $G_1$-invariant alternating bilinear 
form on $V_i$ is zero. 
Since the $V_i$ are mutually non-isomorphic simple $G_1$-modules,
every $G_1$-invariant bilinear pairing
beween $V_i$ and $V_j$ for $i\ne j$ induces
the zero map $V_i \to V_j^* (=V_j)$, and therefore is
zero. Therefore, every $G_1$-invariant
alternating bilinear form on $V=\oplus_{i=1}^n V_i$ is
zero, contradicting that $V$ is symplectic. Thus 
$V_1$ is symplectic.
\end{proof}

We leave the next lemma as an exercise.

\begin{lem}
\label{irredlem}
If $G_0$ is a finite group, $V_0$ is a finite-dimensional
$K$-vector space which is also a 
faithful $K[G_0]$-module,
and $T_0$ is a $G_0$-stable $\O_K$-lattice in $V_0$,
then
$$e_0((x, f),(y,g))= g(x)-f(y)
\text{ for } x, y \in T_0
\text{ and } f, g \in T_0^* := \Hom_{\O_K}(T_0,\O_K)$$ 
defines a  
perfect alternating $G_0$-invariant form on $T_0 \oplus T_0^*$, and 
induces a natural embedding 
$G_0 \hookrightarrow\Aut(T_0 \oplus T_0^*,e_0)$. 
\end{lem}

\begin{lem}
\label{irred}
If the Embedding Theorem is true for all 
irreducible $\rho$ (and $G$) then
it is valid for all $\rho$. 
\end{lem}

\begin{proof}
The $G$-module $V$ splits into a direct sum
of  $G$-modules $V'$ such that
the restriction of $e$ to $V'$ is non-degenerate, and
either $V'$ is simple or $V'=V_0\oplus V_0^*$ where
$V_0$ is simple. 
In the latter case choose a $G$-stable $\O_K$-lattice $T_0$
in $V_0$ and apply Lemma \ref{irredlem}.
\end{proof}

\begin{lem}
\label{Hirred}
Suppose $V$ is simple as both a G-module and an $H$-module
Suppose $T$ is a $G$-stable $\O_K$-lattice in $V$, and
choose $i\in\Z$ so that 
$\pi^ie(T,T)=\O_K$. Then $\pi^ie: T \times T \to \O_K$ 
is a perfect $G$-invariant alternating bilinear form,
and induces an embedding
$$G \hookrightarrow \Aut(T,\pi^ie) \cong \Sp_{2d}(\O_K).$$
\end{lem}

\begin{proof}
Let 
$\bar{e}: T/\pi T \times T/\pi T \to \O_K/\pi\O_K$ be
the non-zero pairing induced by $\pi^ie$. 
Its kernel is an $H$-submodule of $T/\pi T$, 
so is zero by the $H$-simplicity of $T/\pi T$;  
i.e., $\bar{e}$ is non-degenerate. 
By Nakayama's Lemma, $\pi^ie$ is perfect.
\end{proof}

\begin{lem}
\label{extenda}
Suppose $G_0$ is a finite group and $G_1$ is a normal subgroup
of $G_0$. Suppose there exist a free $\O_K$-module 
$T_1$ of rank $2d_1$, 
an alternating perfect form $e_1:T_1\times T_1 \to \O_K$,
and a homomorphism $f:G_1 \to \Aut(T_1,e_1)$
such that whenever $g_1$ is a non-identity element of $G_1$
then there exists $g \in G$ such that $g g_1 g^{-1}\notin\ker(f)$.
Then there exist a free $\O_K$-module $T$ of rank $2d_1[G_0:G_1]$,
an alternating perfect form $e: T \times T \to \O_K$,
and an injective homomorphism $\psi:G_0 \hookrightarrow \Aut(T,e)$.
\end{lem}

\begin{proof}
Let 
$$T=\{u:G_0\to T_1 \mid
u(xs)=s^{-1}u(x) \quad \forall s \in G_1, \forall x \in G_0\},$$
choose a section $p:G_0/G_1 \to G_0$, and let
$$e(u,v)=\sum_{\gamma \in G_0/G_1}e_1(u(p(\gamma)),v(p(\gamma))) 
\quad \text{for } u, v \in T.$$
Note that $e$ is independent of the choice of section $p$. 
Define a homomorphism 
$\psi: G_0 \to \Aut(T,e)$ by 
$\psi(g)(u)(x) =  u(g^{-1}x)$ for 
$g\in G_0$, $u\in T$, $x \in G_0$.
Then the desired conditions are all satisfied.
\end{proof}

\begin{cor}
\label{extendab}
Suppose $G_0$ is a finite group and $G_1$ is a normal subgroup
of $G_0$. Suppose there exists an injective homomorphism
$G_1 \hookrightarrow \Sp_{2d_1}(\O_K)$.
Then  there exists an injective homomorphism
$G_0 \hookrightarrow \Sp_{2d_1[G_0:G_1]}(\O_K)$.
\end{cor}

\begin{lem}
\label{cyclic}
If $\Lambda$ is a finite cyclic group of order $\ell^m$,
then there exists an injective homomorphism
$
\Lambda \times \{\pm 1\} \hookrightarrow
\Sp_{\varphi(\ell^m)}(\Z_\ell)$.
\end{lem}

\begin{proof}
By Lemma 3.7 of \cite{SZcomp} there is an
injective homomorphism 
$\bmu_\ell \times \{\pm 1\} \hookrightarrow
\Sp_{\ell-1}(\Z_\ell)$.
Now apply Corollary \ref{extendab} to
$G_1=\mu_\ell \times \{\pm 1\} \subseteq \Lambda \times \{\pm 1\}=G_0$.
\end{proof}

\section{Proof of the Embedding Theorem}

By Lemmas \ref{irred} and \ref{Hirred},
we may assume from now on that the G-module $V$ is simple and 
the $H$-module $V$ is not simple. 

If (ii) of Proposition \ref{isotype} holds, then 
induct on $\# G$, applying 
Lemma \ref{extenda} to $G_1 \vartriangleleft G_0=G$.

By Proposition \ref{isotype},
we may now assume that $V\cong W^r$ 
for some simple $H$-module $W$, where $r>1$.
It follows that $W$ is self-dual, i.e., there is an $H$-invariant
non-degenerate alternating or symmetric $K$-bilinear form
$f: W \times W \to K$.
We may choose an $H$-stable lattice $T$ in $W$ and (replacing 
$f$ by $\pi^i f$ for suitable $i$,  if necessary) we may assume that
$f: T \times T \to \O_K$
is perfect.
Let
$$w=\dim_K(W).$$

Assume first that $w=1$.
Then
$H \subset \Aut_K(W)=K^*$.
Since $V=W^r$ and $H \subseteq G \subseteq \Sp(V)$,
we have $H \subseteq \{\pm 1\}$, 
and Lemma \ref{cyclic} gives the desired embedding.

Assume from now on that
	$w \ge 2$.
Let 
$$\tau: H \hookrightarrow \Aut(W,f) \subset \Aut(W)$$
be the injective homomorphism defining the $H$-module structure on $W$.
Since $H$ is normal in $G$, the subspace $gW \subset V$ is an
$H$-submodule of $V$ for every $g \in G$. The natural representation
$H \to \Aut(gW)$
is isomorphic to the representation
$$\tau_g: H \to \Aut(W), \quad h\mapsto\tau(g h g^{-1}).$$ 
On the other hand, since $V\cong W^r$ as $H$-modules, therefore 
$gW\cong W$  as $H$-modules. (Indeed, every $H$-submodule
in $V$ is isomorphic to a direct sum of copies of $W$,
and $w=\dim(gW)$.)
Now apply the Extension Theorem and extend $\tau$ 
to a homomorphism
$$\tau_G: G \to \Aut(T,f) \subseteq \Aut(V,f).$$

Since $f$ is perfect, we have $T^*=T$. 
If $f$ is symmetric, then Lemma \ref{irredlem} gives 
a perfect alternating $G$-invariant 
form $e_0$ on $T \oplus T$, and we let
$$\tau_{0}: G \to \Aut(T\oplus T,e_0) \cong \Sp_{2w}(\O_K)$$
be the direct sum of two copies of $\tau_G$.

Suppose $\tau_G$ is injective. 
If $f$ is alternating (resp., symmetric), 
then $\tau_G$ (resp., $\tau_{0}$) gives the desired embedding. 

Now assume $\tau_G$ is not injective. Then $\ker(\tau_G)$
meets $H$ only at the identity, so
$\ker(\tau_G)$ is a normal $\ell$-subgroup of $G$
and therefore is central, and thus contained in 
the Sylow $\ell$-subgroup $L$.
Now retain the notation from the proof of the Extension Theorem.

Suppose $\dim_{E_0}(W) =1$. 
Then $H \subset E_0^*$, and for $h \in H$,
$f(x,y)=f(hx,hy)=f(h^2 x,y) \quad \forall x,y \in W$. 
Thus $h^2=1$, i.e., $H \subseteq \{1,-1\}$, 
and we are done by Lemma~\ref{cyclic}.

Now assume that
	$\dim_{E_0}(W) \ge 2$.
Let 
$$L_0=\{x\in L : \tau_G(x)y=y\tau_G(x) \quad \forall y \in E\}.$$
Then $\ker(\tau_G) \subseteq L_0 \subseteq L$.
By (\ref{idiv}),
\begin{equation}
\label{imiota}
\#\iota(L) \le [E_0:K] \le  w/2.
\end{equation}
Let $G_0= H L_0$, a normal subgroup of $G$. 
Since $V=W^r$, the restriction of $\tau_G$ to $G_0$ 
can be viewed as a homomorphism
$$\tau_G: G_0 \to \Aut_E(W) \subset \Aut_E(V).$$
Then $\rho=\tau_G$ on $H \subset G_0$, and 
$$\rho(z)\rho(h)\rho(z)^{-1}=\tau_G(z)\tau_G(h)\tau_G(z)^{-1}
=\tau_G(z)\rho(h)\tau_G(z)^{-1}$$
for every $z \in L_0, h \in H$. 
Thus  
$$\kappa: L_0 \to \End_H(V)^*=\GL_r(E), \quad 
z\mapsto \tau_G(z)\rho(z)^{-1}$$
is a homomorphism.
Since $\tau_G$ is not injective on $L_0$ but
$\rho$ is, thus
$\kappa$ is injective, since $L_0$ is cyclic. 
Write $\# L_0 = \ell^t$. Then
\begin{equation}
\label{keriota}
r \ge \varphi(\ell^t)=(\ell-1)\ell^{t-1} \ge \ell-1.
\end{equation}
Since
$\# L = \#\ker(\iota)\#\iota(L)=
{\ell}^t \#\iota(L)$,
we have by (\ref{imiota}) and (\ref{keriota}):
$$\varphi({\# L})=\varphi(\ell^t) \#\iota(L) \le 
 r{w}/{2} \le
(r- \frac{\ell -1}{2})w \le 
(r- 2)w
$$
since $\ell \ge 5$. Thus
$$2w +\varphi({\# L}) \le rw =
\dim_K(V) =2d.$$ 
Let 
$$\psi : G/H \cong L \hookrightarrow \Sp_{\varphi(\# L)}(\O_K)$$
be the embedding from Lemma \ref{cyclic}.

If $f$ is alternating (resp., symmetric), 
we are done by taking the direct sum of
$\tau_G$ (resp., $\tau_{0}$) and $\psi$.

\end{document}